\RequirePackage[l2tabu, orthodox]{nag}        
\documentclass[10pt]{article}

\usepackage{graphicx,caption}
\usepackage{booktabs}
\usepackage[T1]{fontenc}
\usepackage{lmodern}
\usepackage[numbers]{natbib}
\usepackage{microtype}
\usepackage{systeme}
\usepackage{ulem}

\usepackage{sistyle}
\SIthousandsep{,}
\usepackage[frenchmath]{mathastext}                    
\usepackage{amsmath}                        
\numberwithin{equation}{section}
\usepackage{amsfonts}
\usepackage{bbold}
\usepackage{commath}
\usepackage{bm}
\usepackage{stackrel}
\usepackage{graphicx}                      
\usepackage[plainpages=false, pdfpagelabels]{hyperref} 
    \hypersetup{
        colorlinks   = true,
        citecolor    = RoyalBlue, 
        linkcolor    = RubineRed, 
        urlcolor     = Turquoise
    }
\usepackage[dvipsnames]{xcolor}
\usepackage[all]{xy}
\usepackage{algorithm}                      
\usepackage{algpseudocode}
\usepackage{amssymb}                        
\usepackage[mathscr]{eucal}                 
\usepackage{dsfont}                                                    
\usepackage[paperwidth=8.5in,paperheight=11.00in,top=1.25in, bottom=1.25in, left=1.00in, right=1.00in]{geometry}
\usepackage{mathtools}                      
\mathtoolsset{showonlyrefs=true}            
\usepackage{fixltx2e,amsmath}               
\MakeRobust{\eqref}
\usepackage{microtype}                                            
\usepackage{amsthm}                         

\usepackage{amsmath}
\usepackage{commath}
\usepackage{algorithm}              
\usepackage{algpseudocode}
\allowdisplaybreaks                         
\theoremstyle{plain}

\numberwithin{theorem}{section}

\theoremstyle{definition}

\newcommand{\X}{\mathcal{X}}

\newcommand{\A}{\mathcal{A}}

\author{
Andrea Angiuli\thanks{Prime Machine Learning Team, Amazon. 320 Westlake Ave N, SEA83, Seattle, WA, 98109 (E-mail: \href{mailto:aangiuli@amazon.com}{aangiuli@amazon.com}). The work presented here does not relate to this author's position at Amazon. }
\and Nils Detering\thanks{Department of Statistics and Applied Probability, South Hall, University of California, Santa Barbara, CA 93106, USA (E-mail: \href{mailto:detering@pstat.ucsb.edu}{detering@pstat.ucsb.edu}). } 
  \and Jean-Pierre Fouque\thanks{Department of Statistics and Applied Probability, South Hall, University of California, Santa Barbara, CA 93106, USA (E-mail: \href{mailto:fouque@pstat.ucsb.edu}{fouque@pstat.ucsb.edu}). Work supported by NSF grants DMS-1814091 and DMS-1953035.} 
  \and Mathieu Lauri\`ere \thanks{Mathematics and Data Science, NYU Shanghai, China (E-mail: 
  \href{mailto:mathieu.lauriere@gmail.com}{mathieu.lauriere@nyu.edu}).}
  \and  Jimin Lin\thanks{Department of Statistics and Applied Probability, South Hall, University of California, Santa Barbara, CA 93106, USA (E-mail: \href{mailto:jiminlin@ucsb.edu}{jiminlin@pstat.ucsb.edu}).} 
  }
\begin{document}

\title{Reinforcement Learning for Intra-and-Inter-Bank Borrowing and Lending Mean Field Control Game}

\maketitle

\begin{abstract}
{We propose a mean field control game model for the intra-and-inter-bank borrowing and lending problem. This framework allows to study the competitive game arising between groups of collaborative banks. The solution is provided in terms of an asymptotic Nash equilibrium between the groups in the infinite horizon. A three-timescale reinforcement learning algorithm is applied to learn the optimal borrowing and lending strategy in a data driven way when the model is unknown. An empirical numerical analysis shows the importance of the three-timescale, the impact of the exploration strategy when the model is unknown, and the convergence of the algorithm.}
\end{abstract}

\section{Introduction}
Many problems in finance involve a large number of strategic agents. A typical example is how traders interact in a common market through the price of some assets. At a larger scale, another example is how banks interact through the money they borrow from or lend to each other or to a central bank. When the agents are competing, one can represent the problem as a game and look for a Nash equilibrium. On the other hand, when the agents are cooperating, one can look for a social optimum. The problems have different solutions, and the non-cooperative equilibrium has a higher average cost per player, which is interpreted as a lack of efficiency. This leads to the notion of price of anarchy~\cite{koutsoupias1999worst}. 

When the number of agents is large, studying every pairwise interactions becomes intractable. To simplify the analysis, a mean field approximation can be used, assuming that the population is homogeneous and the interactions are symmetric. This idea led to the notion of mean field games (MFGs) and mean field control (MFC) problems (also known as McKean Vlasov control problems) depending on whether the agents are competitive or cooperative, \cite{MR2295621,MR2346927,bensoussan2013mean,carmona2018probabilisticI-II}. A related notion is the concept of mean field type game (MFTG), in which a finite number of players compete and each player's problem is of MFC type~\cite{tembine}. MFG, MFC and MFTG have found applications in energy production and management~\cite{gueant2011meanappli,alasseur2020extended,chan2015bertrand}, crowd trading~\cite{cardaliaguetlehalle2018mean}, systemic risk~\cite{CarmonaFouqueSun}, to cite just a few examples. See e.g.~\citep{carmona2020applicationsfinancial} for a recent survey of applications to finance and economics. 

In the past few years, the question of \textit{learning} solutions to MFG and MFC problems using model-free methods based on reinforcement learning (RL) has gained momentum. Many of these methods rely on updating a value function and a distribution. In particular, stationary MFG solutions have been approximated in \cite{guo2019learning} using fixed point iterations and Q-learning, and in \cite{elie2020convergence} using fictitious play and deep RL. Two-timescale analysis to learn MFG solutions have been used in~\cite{mguni2018decentralised,SubramanianMahajan-2018-RLstatioMFG}.  Recently, a two timescale algorithm has been introduced in~\cite{AFL2021} to solve MFG or MFC depending on the choice of learning rates for the distribution and the value function. This allows to have a unified point of view on these two types of problems and a common RL method. In \cite{ADFLL-MFCGpaper}, the approach has been extended to mean field control games (MFCG) using a three-timescale RL algorithm. It was developed in a finite-horizon setting for extended MFCGs arising naturally in the context of the trader's liquidation problem between competitive groups of collaborative traders who share the inventory cost of their group. 

In the present paper, our main contributions are threefold. First, in Section~\ref{sec:bank-model}, we introduce a model of intra-and-inter-bank borrowing and lending, which can be viewed as an extension of the model studied in~\cite{CarmonaFouqueSun} where there are local coalitions inside each bank. Second, in Section~\ref{sec:3timescales}, we apply a three-timescale RL algorithm to solve this class of infinite horizon problems. Last, in Section~\ref{sec:numerics}, we show numerical results that illustrate the performance of our method on the model of intra-and-inter-bank borrowing and lending.

\section{Intra-and-Inter-Bank Borrowing and Lending Problem}
\label{sec:bank-model}
A model of inter-bank borrowing and lending has been introduced in \cite{CarmonaFouqueSun} as a linear-quadratic stochastic differential game between banks which control their drifts and minimize a quadratic cost with incentive to stay close to the average capitalization of the system. The model has been studied as a finite-player game in finite horizon. Open-loop and closed-loop Nash equilibria have been computed using Forward-Backward Stochastic Differential Equations (FBSDE) and Hamilton-Jacobi-Bellman (HJB) partial differential equations. In this model the central bank acts as the clearing house. Systemic risk has then been considered as a large deviation event. In addition to the finite player game a Mean Field Game (MFG) limit has been discussed as well. 
In the present paper we first propose an extension of the aforementioned model where the competitive banks are made of collaborative branches leading to a Mean Field Control Game (MFCG) model, and, second, we use a \textit{three-timescales reinforcement learning algorithm} to solve this problem when the structures of the dynamics and of the cost are unknown to the agents. This represents a natural and interesting development of the \textit{two-timescales reinforcement learning algorithm} introduced in \cite{AFL2021} to solve MFG or Mean Field Control (MFC) problems. The following model of intra-and-inter-bank borrowing and lending provides a benchmark for our algorithm, which can be applied to a wide range of models.

\subsection{System of \texorpdfstring{$M$}{M} Banks with \texorpdfstring{$N$}{N} Branches }

In the model considered below, we consider $M \in \mathbb{N}$ bank groups. Each bank has $N \in \mathbb{N}$ local branches and is involved in both inter- and intra- bank borrowing and lending activity. Let tuple $(m,n)$ for $m \in \{1, \dots, M\}$ and $n \in \{1, \dots, N\}$ index the $n$-th branch of the $m$-th bank. The one-dimensional diffusion process $(X_t^{m,n})_{t \in [0, \infty)}$ stands for the \textit{log-monetary reserve} of the branch $(m, n)$ over an infinite time horizon, whose dynamics has the following form:
\begin{align}
d X_t^{m,n}
    = \left[\kappa \left( \frac{1}{N}\sum_{j=1}^{N} X_t^{m,j} - X_t^{m,n} \right) + \alpha_t^{m}(X_t^{m,n})\right]dt + \sigma d W_t^{m,n},
\end{align}
with $X_0^{m,n} \sim \mu_0$. The first term in the drift
$$
\kappa \left( \frac{1}{N}\sum_{j=1}^{N} X_t^{m,j} - X_t^{m,n} \right)
=\frac{\kappa}{N}
\sum_{j=1}^{N}\left(  X_t^{m,j} - X_t^{m,n} \right), \quad \kappa\geq 0,
$$
represents borrowing and lending activity between branch $(m,n)$ with the other branches of the same bank. As can be seen on the right-hand side, branches with more liquidity lend to branches with less liquidity at a rate $\kappa$, normalized by the number $N$ of branches. The left-hand side can be interpreted as mean-reversion to the average liquidity reserve of the branches of that bank. To some degree (depending on $\kappa$) this mean reversion will be facilitated by the branches at no cost because branches that are well equipped with liquidy have an interest in investing their excess liquidity and branches with too little liquidity have an interest in borrowing.
 In addition to this mean-reversion behavior, local branch $(m,n)$ has the possibility to borrow and lend from a central bank. This borrowing happens at a rate that depends on the liquidity reserve of $(m,n)$ but needs to comply with the (time-dependent) feedback-form policy of bank group $m$, which is reflected in the \textit{control} term $\alpha^m: \mathbb{R} \rightarrow \mathbb{R}$. The entire system is driven by $M \times N$ independent standard Brownian motions $(W_t^{m,n})_{t,m,n}$. For simplicity, we assume the same constant diffusion rate $\sigma > 0$.
Bank group $m$ designs its policy of control $\alpha_t^m$ of the borrowing and lending rate for all of its branches at time $t$ in order to minimize the \textit{group objective function}
\begin{align}
\label{total-cost-group-m}
J(\alpha^m; \alpha^{-m})
    = \frac{1}{N} \sum_{n=1}^{N} \mathbb{E}\left\{ \int_{0}^{\infty}e^{-\beta t} f\left(X_t^{m,n}, \alpha_t^m(X_t^{m,n}), \mu_t, \widetilde\mu^m_t\right) dt\right\},
\end{align}
where $\beta>0$ denotes the time discount rate, and the interaction with the other banks is through the \textit{global empirical distribution} $\mu_t = \frac{1}{MN} \sum_{i=}^{M}\sum_{j=1}^{N} \delta_{X_t^{i,j}}$ of reserves of the entire system across all branches and all banks, while the interaction within the branches of bank $m$ is through the \textit{local empirical distribution} $\widetilde\mu^m_t = \frac{1}{N} \sum_{j=1}^{N} \delta_{X_t^{m,j}}$. We denote by $\alpha^{-m}$ the control profile for all bank groups except $m$, i.e., $\alpha^{-m} = (\alpha^1,\dots,\alpha^{m-1},\alpha^{m+1},\dots,\alpha^M)$. Here, we consider a quadratic  \textit{running cost} function given by
\begin{align}\label{runningcost}
f(x, \alpha, \mu, \widetilde{\mu})
    = \frac{1}{2}\alpha^2 + c_1(x-c_2 \overline{\mu})^2 + \widetilde{c}_1(x-\widetilde{c}_2 \overline{\widetilde{\mu}})^2 + \widetilde{c}_3(\overline{\widetilde{\mu}} - \widetilde{c})^2,
\end{align}
which depends on the global and local empirical distributions $\mu, \widetilde{\mu} \in \mathcal{P}(\mathbb{R})$ only through their first moments, denoted respectively by $\overline{\mu}, \overline{\widetilde{\mu}}$. So in~\eqref{total-cost-group-m}, the cost at time $t$ depends only on the global and local empirical means 
$\overline{\mu}_t = \frac{1}{MN} \sum_{i=}^{M}\sum_{j=1}^{N} X_t^{i,j}$ and 
$\overline{\widetilde\mu}^m_t = \frac{1}{N} \sum_{j=1}^{N} X_t^{m,j}$.
Here, $c_1, c_2, \widetilde{c}_1, \widetilde{c}_2, \widetilde{c}_3, \widetilde{c} \in \mathbb{R}$ are some constants. The running cost is interpreted as follows. The first term represents the quadratic cost of control on borrowing and lending rate. The second and third term shows the bank's intention to keep the reserve of its branch close to both global average reserve $\overline{\mu}$ and local average reserve $\overline{\widetilde{\mu}}$ to some extend quantified by $c_1, c_2, \widetilde{c}_1, \widetilde{c}_2$. Meanwhile, the bank prefers its local average centering around a target level $\widetilde{c}$. 

The above system of $M$ banks constitutes a competitive game between the $M$ banks, while it is a collaborative (distributed) game within each bank group. We are looking for a closed-loop Nash equilibrium between the banks. This kind of mixed competitive-collaborative game is described in \cite{ADFLL-MFCGpaper} in the context of finite horizon extended games applied to the liquidation trader's problem. Mathematically, the problem is defined as follows: find a control profile $(\hat\alpha^m)_{m}$ such that for every $m=1,\dots,M$, $\hat\alpha^m$ minimizes:
$$
    \alpha^m \mapsto J(\alpha^m; \hat\alpha^{-m}).
$$
We are interested in the \textit{Mean Field Control Game Limit} when both $M$ and $N$ go to infinity and its solution using reinforcement learning as presented in Section \ref{sec:3timescales}.

\subsection{Mean Field Control Game Limit}

Associated to the finite-player game introduced above, we associate the Mean Field Control Game (MFCG) obtained in the asymptotic 
limit where both $M$ and $N$ go to $\infty$. The problem is to find a pair $(\hat\alpha,\hat\mu)$ such that the following two conditions are satisfied:
\begin{enumerate} 
\item A representative bank confronted with a fixed flow of probability distributions $\hat\mu:=(\hat\mu_t)$ solves the McKean-Vlasov (MKV) control problem of finding a minimizer $\hat{\alpha}$ for
\begin{align} \label{eq:mfcg_obj}
\alpha \mapsto J(\alpha; \hat\mu)
    = \mathbb{E}\left\{ \int_{0}^{\infty}e^{-\beta t} f\left(X_t, \alpha_t(X_t), \hat\mu_t, \mathcal{L}(X_t) \right) dt \right\},
\end{align}
subject to
\begin{align} \label{eq:mfcg_dx}
dX_t = \left[\kappa(\mathbb{E}(X_t) - X_t) + \alpha_t(X_t) \right]dt + \sigma dW_t, \quad X_0 \sim \mu_0.
\end{align}
\item The law of the state $X_t$ controlled by $\hat\alpha$ satisfies the fixed point condition
\begin{align} \label{eq:mfcg_fp}
 \mathcal{L}(X_t) = \hat\mu_t, \quad t\in[0,\infty).
\end{align}
\end{enumerate}
The justification of such a limit is treated mathematically in the forthcoming paper \cite{ADFLproof}. See also the Appendix in \cite{ADFLL-MFCGpaper} for a formal justification in the case of the linear-quadratic trader's liquidation problem.

\subsubsection{Value Function and HJB Equation}

Since we are looking for an equilibrium among Markovian feedback strategies, we solve the MFCG system \eqref{eq:mfcg_obj}-\eqref{eq:mfcg_fp} through the  Hamilton-Jacobi-Bellman (HJB) equation approach. 
Following the computation detailed in the Appendix A of \cite{AFL2021}, we first solve the finite horizon problem with zero terminal condition when the global distribution flow is given by $(\mu_t)_{t \in [0,T]}$:
\begin{align}
V^T(t,x)=\inf_{\alpha}\mathbb{E}\left\{ \int_{t}^{T}e^{-\beta s} f\left(X_s, \alpha_s(X_s), \mu_s, \mathcal{L}(X_s) \right) ds \right\},
\end{align}
subject to: 
\begin{align*}
    dX_s = \left[\kappa(\mathbb{E}(X_s) - X_s) + \alpha_s(X_s) \right]ds + \sigma dW_s, \quad X_t =x
\end{align*}
and with the fixed point condition \eqref{eq:mfcg_fp} over $[t,T]$.
Denoting by $\mathcal{A}$ the infinitesimal generator of $X$, the Hamiltonian is given by
\begin{align}
H(t, x, \hat{\alpha}(t,x), \mu_t,  \widetilde{\mu}_t)
    &= \inf_{\alpha}\left\{ \mathcal{A}V^T(t,x) + f(x, \alpha, \mu_t, \widetilde{\mu}_t) \right\},
\end{align}
which attains its minimum at $\hat\alpha(t,x) =- \partial_x V^T(t,x)$ in our case where $f$ is given by \eqref{runningcost}, and the dynamics of $X$ by \eqref{eq:mfcg_dx}.
The HJB equation with MKV dynamic reads (see e.g.~\cite[Section 4.1]{bensoussan2013mean})
\begin{align} \label{eq:hjb}
    \partial_t V^T(t,x) -\beta V^T(t,x) + H(t, x, \hat{\alpha}(t,x), \mu_t,  \widetilde{\mu_t}) + \int_{\mathbb{R}}\frac{\partial H}{\partial \widetilde{\mu}_t}(t, \xi, \hat{\alpha}(t,\xi), \mu_t,  \widetilde{\mu}_t)(x) d \widetilde{\mu}_t(\xi)= 0,
\end{align}
with $V^T(T,x)=0$.

We compute
\begin{align*}
    H(t, x,\hat\alpha(t,x),\mu_t,\widetilde\mu_t)
    &= -\frac{1}{2}(\partial_x V^T(t,x))^2+\frac{1}{2}\sigma^2 \partial_{xx}V^T(t,x) +\kappa(\overline{\widetilde{\mu}}_t-x)V^T_x(t,x) \\
    &\quad +c_1(x-c_2 \overline{\mu}_t)^2 + \widetilde{c}_1(x-\widetilde{c}_2 \overline{\widetilde{\mu}}_t)^2 + \widetilde{c}_3(\overline{\widetilde{\mu}}_t - \widetilde{c})^2,
\end{align*}
and
\begin{align*}
       \int_{\mathbb{R}}\frac{\partial H}{\partial \widetilde{\mu}_t}(t,\xi, \hat{\alpha}(t,\xi), \mu_t, \widetilde{\mu}_t)(x) d \widetilde{\mu}_t(\xi) 
       &=
     -2\widetilde{c}_1\widetilde{c}_2(1-\widetilde{c}_2)\overline{\widetilde{\mu}}_t x + 2\widetilde{c}_3(\overline{\widetilde{\mu}}_t - \widetilde{c})x.
\end{align*}
We then formulate the following ansatz for the value function
\begin{align}
V^T(t,x) = \Gamma^T_2(t) x^2 + \Gamma_1^T(t) x + \Gamma^T_0(t),
\end{align}
with the zero terminal conditions $\Gamma^T_2(T)=\Gamma^T_1(T)=\Gamma^T_0(T)=0$. We have $\hat\alpha(t,x)=-2\Gamma^T_2(t)x-\Gamma_1^T(t)$. Plugging the ansatz and its partial derivatives
into \eqref{eq:hjb} and identifying the coefficients of powers of $x$ 
leads to a system of ODEs for $\Gamma^T_1,\Gamma^T_2,\Gamma^T_0$ with zero terminal conditions. This system is complemented with the forward equation
\begin{align}\label{eq:mu}
    d\overline\mu_t
    = \mathbb{E}\left(\hat\alpha(t,X_t)\right)dt=-\left[2\Gamma^T_2(t)\overline\mu_t+\Gamma_1^T(t)\right]dt, \quad \overline\mu_0=x,
\end{align}
obtained by taking expectation in \eqref{eq:mfcg_dx} and using the expression of the control $\hat\alpha$.
The ODE system for $(\Gamma^T_2(t),\Gamma^T_1(t),\Gamma^T_0(t),\overline\mu_t)_{t \in [0,T]}$ is a two-point boundary value problem which can be solved explicitly as in in the Appendix A of \cite{AFL2021}. 
\subsubsection{Explicit Formulas}\label{subsec: benchmark}
The solution to our infinite horizon problem is obtained by taking the limit $T\to \infty$. Furthermore, since we are interested in the \textit{asymptotic solution}, or equivalently the \textit{stationary solution}, we take the limit $t\to\infty$ to obtain that the limiting value function
\begin{align*}
V(x)&= \Gamma_2 x^2 + \Gamma_1 x + \Gamma_0,
\end{align*}
where $\Gamma_0, \Gamma_1, \Gamma_2$ are constants, 
must satisfy \eqref{eq:hjb} with $\partial_t V^T=0$, no terminal condition at $T=+\infty$, and  $\overline{\hat\mu}_t=\overline{\hat\mu}$ being the stationary point of \eqref{eq:mu} satisfying $2\Gamma_2\overline{\hat\mu}+\Gamma_1=0$. We deduce the formulas:
\begin{align*}
\hat\alpha(x)&=-2\Gamma_2 x - \Gamma_1 ,\\
\Gamma_2&=\frac{-(\beta + 2\kappa) + \sqrt{(\beta + 2\kappa)^2 + 8(c_1 + \widetilde{c}_1)}}{4},\\
\Gamma_1&= \frac{2\widetilde{c}_3(\overline{{\mu}} - \widetilde{c}) - 2\widetilde{c}_1\widetilde{c}_2(2-\widetilde{c}_2)\overline{{\mu}} - 2 c_1c_2\overline{\mu}}{\beta + \kappa + 2\Gamma_2},\\
\Gamma_0&= \frac{ -\kappa\overline{{\mu}} - \frac{1}{2}\Gamma_1^2 + \sigma^2\Gamma_2 + c_1c_2^2 \overline{\mu} + \widetilde{c}_1\widetilde{c}_2^2\overline{{\mu}} + \widetilde{c}_3(\overline{{\mu}} - \widetilde{c})^2}{\beta},\\
\overline\mu&= -\frac{\Gamma_1}{2\Gamma_2}
    = \frac{\widetilde{c}_3 \widetilde{c}}{c_1(1-c_2) + \widetilde{c}_1(1-\widetilde{c}_2)^2 + \widetilde{c}_3 - \kappa \Gamma_2}.
\end{align*}
Note that at Nash equilibrium and asymptotically when time is large, $X_t$ behaves like an OU process with a rate of mean-reversion $\kappa+2\Gamma_2$ and diffusion $\sigma^2$. Therefore, the equilibrium asymptotic distribution is $\mu=\mathcal{N}\left(\overline\mu, \frac{\sigma^2}{2\kappa+4\Gamma_2}\right)$.

\section{Three-Timescale \texorpdfstring{$Q$}{Q}-Learning Algorithm}\label{sec:3timescales}

\subsection{Discrete time formulation and \texorpdfstring{$Q$}{Q}-learning}
We now describe our algorithm to learn the solution to the mixed Control Game problem  \eqref{eq:mfcg_obj}-\eqref{eq:mfcg_fp}. Since the algorithm itself is only a minor modification of the algorithm used in \cite{ADFLL-MFCGpaper}, we keep this paragraph brief. The algorithm rests on the concept of $\mathcal{Q}$ learning, a well established method to solve Markov Decision problems. We first discretize the time interval $[0,\infty]$ into an equally spaced grid $0=t_0 < t_1 < \dots $ and assume for notational simplicity that $t_i=i$. We then recast the problem \eqref{eq:mfcg_obj}-\eqref{eq:mfcg_fp} into a discrete time mean field control game problem given by:
\begin{enumerate}
\item Given $\{\mu_n\}_{n\in \mathbb{N}}$, find a minimizer $\hat{\alpha }$ for
\begin{align}\label{disc:opt}
    J(\alpha; \mu ) & =\mathbb{E}\left[\sum_{n=0}^{\infty} e^{-\beta n} f\left(X_{n},\alpha_n(X_{n}),\mu_n,\mathcal{L} (X_{n})\right)\right],
\end{align}
subject to 
\begin{align*}
    \mathbb{P}\big(X_{n+1} = x' | X_{n} = x, && \alpha_n (X_{n}) = a, &&\mu_n = \mu,  && \mathcal{L} (X_{n}) = \widetilde{\mu} \big)
    = p(x'|x,a,\mu, \widetilde{\mu}),
\end{align*}
where the transition kernel $p:\X \times \A \times \Delta^{ | \X |} \times \Delta^{ | \X |} \rightarrow \Delta^{ | \X |}$ arises from a discrete counterpart to  \eqref{eq:mfcg_dx}. \item The law of the state $X_n$ matches the fixed point condition
\begin{align} \label{disc:fp}
 \mathcal{L}(X_n) = \mu_n, \quad n\in \mathbb{N}.
\end{align}
\end{enumerate}

In order to solve this discrete time problem we discretize the state space into $\X = \{x_0, \dots, x_{\abs{\X}-1}\}$, and action space into $\A = \{a_0, \dots, a_{\abs{\A}-1}\}$ respectively. 

Our reinforcement learning algorithm to solve the discrete time and discrete state problem \eqref{disc:opt} and \eqref{disc:fp} follows \cite{ADFLL-MFCGpaper}. The algorithm is based on well established ideas from $Q$-learning. The algorithm is model agnostic which means that no information is needed about the model that generates the data. In the control part of our problem \eqref{disc:opt}, the local distribution $\mathcal{L}(X_n)$ depends on the control that is chosen. For this reason it can not simply be treated as an additional parameter but the $Q$ learning has to be adapted slightly. For an admissible control $\alpha : \X \rightarrow \A$, we define the new control $\alpha_{x,a}$ that deviates from $\alpha$ only at the state $x$ where it takes the value $a$:
\begin{equation}
\alpha_{x,a} (x')=\begin{cases}
      a & \text{if   } x' = x \\
      \alpha(x)  & \text{otherwise}.
\end{cases}\label{mod:alpha}
\end{equation}
Given a fixed global measure $\mu$ and strategy $\alpha$, the $Q$-function for our problem is then given by:
\begin{align}
    Q^{\alpha}_{\mu} (x,a) = f(x,a,\mu,\mu^{\alpha_{x,a}}) + \mathbb{E}\left[\sum_{n=1}^{\infty} e^{-\beta n}f(X_{n},\alpha(X_{n}),\mu,\mu^\alpha)\lvert X_{0}=x,A_{0}=a\right]. 
\end{align}
One can then consider the optimal cost function 
$$
    Q_{\mu}^{*} (x,a) := \min_{\alpha} Q_{\mu}^{\alpha}(x,a),
$$ 
which, conditioned on being in state $x$ and choosing action $a$ at time $0$, minimizes the cost over all strategies $\alpha$ chosen in all steps to follow.
From the function $Q_{\mu}^{*}$ one obtains the optimal control $\alpha^* (x)=\arg \min_a Q_{\mu}^{*} (x,a)$. In Section~\ref{sec: action_exploration} we will see that actually a randomized counterpart of $\alpha^*$ should be chosen to ensure a wide enough exploration range of the possible actions. We stress that the minimizing strategy usually depends on the global measure $\mu$. For fixed $\mu$, it follows from \cite{AFL2022}, as the measure $\mu$ is fixed and does not depend on $\alpha$, that the function $Q_{\mu}^{*}$ follows a Bellman equation given by: 
\[
    Q_{\mu}^{*}(x,a)=f(x,a,\mu,\mu_{x,a}^{*})+\gamma\sum_{x'}p(x'\lvert x,a,\mu,\mu_{x,a}^{*})\min_{a'}Q_{\mu}^{*}(x',a').
\]
The measure $\mu_{x,a}^{*} = \lim_{n\rightarrow\infty}\mathcal{L} (X_{n}^{\alpha^*_{x,a},\mu})$ corresponds to the strategy $\alpha^*_{x,a}$ as derived from $\alpha^*$ by changing the action in state $x$ to $a$, see \eqref{mod:alpha}.

\subsection{Three-Timescale Updating Rates}
Our algorithm to approximate the $Q$ function, optimal policy and the equilibrium distribution mimics the idea of nested optimization. For a given global distribution, the $Q$-function that describes the optimal action has to be found, and based on this, the local distribution. This idea of nested simulation leads to a Three-Timescale approach which is sketched in the following. With updating rates $\rho_{k}^{\mu}$ for the global distribution, $\rho_{k}^{Q}$ for the $Q$ table, and $\rho_{k}^{\mu^{\alpha}}$ for the local distribution, where we assume $\rho_{k}^{\mu} < \rho_{k}^{Q} < \rho_{k}^{\mu^{\alpha}}$, the updates that can be derived from the Bellman equation are described by
\begin{align}\label{update:system}
\begin{cases}
\mu_{k+1}
    =\mu_{k}+\rho_{k}^{\mu}\mathcal{P}(Q_{k},\mu_{k}),\\
\mu^{\alpha}_{k+1}
    =\mu^{\alpha}_{k}+\rho_{k}^{\mu^{\alpha}}\mathcal{P}(Q_{k},\mu^{\alpha}_{k}), \\
Q_{k+1}
    =Q_{k}+\rho_{k}^{Q}\mathcal{T}(Q_{k},\mu_{k},\mu^{\alpha}_{k}),
\end{cases}
\end{align}
with
\begin{align}
\begin{cases}
\mathcal{P}(Q,\nu)(x)
    =(\nu P^{Q,\mu,\mu^{\alpha}})(x)-\nu(x), \\
\mathcal{T}(Q,\mu,\mu^{\alpha})(x,a)
    =f(x,a,\mu,\mu^{\alpha}) +\gamma\sum_{x'}p(x'\lvert x,a,\mu,\mu^{\alpha})\min_{a'}Q(x',a')-Q(x,a)\\
P^{Q,\mu,\mu^{\alpha}}(x,x')
    =p(x'\lvert x,\arg\min_{a}Q(x,a),\mu,\mu^{\alpha}), \\
(\nu P^{Q,\mu,\mu^{\alpha}})(x)
    =\sum_{x_{0}}\nu(x_{0})P^{Q,\mu,\mu^{\alpha}}(x_{0},x).
\end{cases}
\end{align}
Note that in our model agnostic approach, the transition probabilities $p$ need to be estimated from the data. As samples from the state and the rewards are obtained incrementally, we update these estimates with Robbins–Monro rates. We refer the reader to \cite{ADFLL-MFCGpaper} for more details. 

\subsection{Action Exploration} \label{sec: action_exploration}
An efficient algorithm is designed to well balance the tendencies between exploring a range of policies and staying in the current best choice, i.e. exploration and exploitation. An over-exploring algorithm is less likely to converge to the optimal policy while the over-exploiting one will possibly be stuck in a local optimal, which is the well known exploration-exploitation dilemma \cite{kaelbling1996reinforcement}. As other reinforcement learning algorithms, our three-time scale $Q$-learning algorithm is confronted with this dilemma. Therefore, we shall develop methods to balance the exploration-exploitation trade-off.

Over the recent decades, various action exploration techniques have been developed to overcome the exploration and exploitation dilemma. Those can roughly be distinguished into two categories: \textit{undirected} and \textit{directed} \cite{thrun1992efficient}. Undirected exploration takes actions based on some probability distribution and does not account for the learning progress itself. Widely applied undirected methods include $\epsilon$-greedy, Boltzmann, and Max-Boltzmann \cite{wiering1999explorations}. In contrast, directed exploration adapts the action preference by the learning progress, such as the number of times of a state-action pair being visited (counter-based), the environment with large errors from previous exploration (error-based), states not being visited recently (recency-based). 

Depending on specific learning tasks, sophisticated directed exploration might require more efforts to calibrate but does not necessarily outperform simple undirected heuristics \cite{wiering1998efficient, kuleshov2014algorithms}. Therefore, for the new three-timescale algorithm that has not been comprehensively tested, we shall first focus on the undirected methods, with preference for its generality and simplicity. It can then serve as a benchmark for the application of more complicated directed exploration methods. In particular we consider the following three undirected exploration methods:
\begin{enumerate}
\item \textit{$\epsilon$-greedy}.
    \begin{align} \label{eq:greedy}
        \pi_t^{\epsilon}(x)
        = \begin{cases}
        a \sim \text{Unif}(\mathcal{A}), & \text{ w.p. } \epsilon,\\
        \arg\max_{a \in \mathcal{A}}Q_t(x,a), & \text{ w.p. } 1 - \epsilon.
        \end{cases}
    \end{align}
    Parameter $\epsilon$ is the \textit{exploration rate}.
\item \textit{Boltzmann exploration}.
    \begin{align} \label{eq:boltzmann}
        \pi_t^{Boltz}(x, a) \sim Boltz(Q_t(x,a);Q_t(x,\cdot), \tau)
    \end{align}
    with $Boltz(x; X, \tau):=\frac{e^{-x/\tau}}{\sum_{x'\in X} e^{-x / \tau}}$ known as the Boltzmann distribution. Parameter $\tau$ is referred as the \textit{temperature}.
\item \textit{Max-Boltzmann} combines the $\epsilon$-greedy with Boltzmann exploration by replacing $\text{Unif}(\mathcal{A})$ in \eqref{eq:greedy} by $Boltz$ distribution in \eqref{eq:boltzmann},
\end{enumerate}
where the exploration propensity of the algorithm is controlled by the exploration rate $\epsilon$ or constant temperature $\tau$. To search for the appropriate exploration heuristic, for each of the three heuristics, we consider the following three configurations: (1) constant rate; (2) linearly decaying rate w.r.t episode; and (3) exponentially decaying rate w.r.t episode, which will be specified in Section~\ref{sec:numerics}.

\subsection{Algorithm}
The Algorithm~\ref{alg:rl_infinite} applied to learn the asymptotic solution discussed in section \ref{subsec: benchmark} is the three-timescale mean field Q-learning algorithm (U3-MF-QL) presented in \cite{ADFLL-MFCGpaper}. By interacting with the environment in a trial and error fashion, we are able to learn the optimal $Q$ table, together with the local and global distribution at equilibrium. As discussed in the previous section, the learning rates assume a core role and they are defined as 
\begin{align}\label{eq:rates}
\rho_{x,a, n, k}^{Q}:=\frac{1}{(1+ \#\abs{(x, a, k, n)})^{\omega^Q}}, &&
\rho_{k}^{\nu}:=\frac{1}{(1+k)^{\omega^{\nu}}},
\end{align}
where $\nu$ is replaced by $\mu$ and $\widetilde\mu$ for the local and global distribution respectively, and $\#\abs{(x, a, k, n)}$ counts the visits of the pair $(x,a)$ up to the episode $k$ and time $n$. The triplet $(\omega^{Q}, \omega^{\mu}, \omega^{\widetilde{\mu}})$ should be chosen such that $\omega^{\mu} > \omega^{Q} >\omega^{\widetilde{\mu}}$, so that $\rho^{\nu}_k < \rho^Q_k < \rho^{\widetilde{\nu}}_k$, and it should satisfy $\omega^{Q} \in (0.5, 1)$. 
\begin{algorithm}[H]
\caption{Three-Timescales Mean Field Q-Learning - Infinite Horizon}
\begin{algorithmic}[1]
\Require{\\
    \qquad T: number of time steps in a learning episode, \\
    \qquad Truncated state space: $\X = \{x_0, \dots, x_{\abs{\X}-1}\}$, \\
    \qquad Truncated action space: $\A = \{a_0, \dots, a_{\abs{\A}-1}\}$, \\
    \qquad Initial distribution of the representative player: $\mu_0$, \\
    \qquad Exploration rule s.t. $\pi^v \in \Delta^{|\A|}$ for any $|\A|-$dim vector $v$, \\
    \qquad Break rule tolerances: $tol_Q$, $tol_{\mu}$, $tol_{\widetilde{\mu}}$.
}
\State{\textbf{Initialization:}}
    \State{\qquad $Q^0(x, a)=0$ for all $(x,a) \in \X \times \A$,}
    \State{\qquad $\mu^0_{n} = \frac{1}{\abs{\X}}J_{\abs{\X}}$ and $\widetilde{\mu}^0_{n} = \frac{1}{\abs{\X}}J_{\abs{\X}}$ for $n = 0,\dots, T$,}
    \State{\qquad where $J_{m}$ is an $m$-dimensional unit vector.}
\For{each episode $k=1, 2, \dots$}
\State{\textbf{Set }$Q^k \equiv Q^{k-1}$}
    \State{\textbf{Observe initial state:} $X_{0}^k \sim \mu_T^{k-1}$.}
    \For{$n = 0, \dots, {T}$}
        \State{\textbf{Choose action:}}
            \State{choose $A_{n}^k$ using the exploration policy $\pi^{Q^{k}(X_{n}^k, \cdot)}$.}
        \State{\textbf{Update distributions:}}
            \State{\qquad $\mu^{k}_{n} = \mu^{k-1}_{n} + \rho_{k}^{\mu}(\boldsymbol{\delta}(X_{n}^k) - \mu^{k-1}_{n})$,}
            \State{\qquad $\widetilde{\mu}^{k}_{n} = \widetilde{\mu}^{k-1}_{n} + \rho_{k}^{\widetilde{\mu}}(\boldsymbol{\delta}(X_{n}^k) - \widetilde{\mu}^{k-1}_{n})$,}
            \State{\qquad where $\boldsymbol{\delta}(X_{n}^k)=\left(\mathbf{1}_{x}(X_{n}^k)\right)_{x \in \X}$.}
        \State{\textbf{Observe next state:}}
            \State{\qquad observe $X_{{n+1}}^k$ from the environment.}
        \State{\textbf{Observe cost:}}
            \State{\qquad observe $f_{n} = f(X_{n}^k, A_{n}^k, \mu^{k}_{n}, \widetilde{\mu}^{k}_{n})$.}
        \State{\textbf{Update $Q$ table:}}
            \State{\qquad $\begin{aligned} Q^{k}(x,a) =& \, Q^{k}(x,a) +  \mathbf{1}_{x,a}(X_{n}^k, A_{n}^k)\rho^{Q}_{x,a, n, k}\\
            &(f_{n} +
            \beta\min_{a' \in \A} Q^{k}(X_{{n+1}}^k, a')-Q^{k}(x,a)),\end{aligned}$}
         \State{\qquad where $\beta$ is the discount parameter.}
    \EndFor
    \If{$\begin{cases}
    \delta(\mu^{k},\mu^{k-1}) \le tol_{\mu},\\
    \delta(\widetilde{\mu}^{k},\widetilde{\mu}^{k-1}) \le tol_{\widetilde{\mu}},\\
    \Vert Q^k - Q^{k-1} \Vert_{1,1} \leq tol_{Q}, 
    \end{cases}$}{ break}
    \EndIf
\EndFor
\end{algorithmic}
\label{alg:rl_infinite}
\end{algorithm}

\section{Numerical Results}
\label{sec:numerics}
For the MFCG problem setting, we choose $(c_1, c_2, \widetilde{c}_1, \widetilde{c}_2, \widetilde{c}_3, \widetilde{c}) = (1.5, 0.75, 2.5, 0.5, 4, 2)$ and discount rate $\beta=1$ for the running cost $f$; $(\kappa, \sigma) = (1, 2)$ for the dynamic of state $dX$. We truncate the infinite time horizon by $[0, T]$ with $T=20$ and discretize it by steps of size $\delta t = 1/16$. The state and action spaces are trimmed into $\X = \{x_0=-1.5, x_1=-1.5+\delta x, \dots, x_{\abs{\X}-1}=4.5\}$ and $\A = \{a_0=-6, a_1=-6 + \delta a, \dots, a_{\abs{\A}-1}=6\}$ by $\delta x = \delta a = \sqrt{\delta t} = 1/4$. For the reinforcement learning setup, we take $K=50,000$ episodes and consider the specifications for the action exploration in Table~\ref{tab:act_expl}. The initial exploration rate is set small for the constant $\epsilon$-greedy action explorer, mildly greater for the linearly decaying rate, and large for the exponentially decaying rate. The initial temperature for Boltzmann explorers are the same. The Max-Boltzmann explorers takes in a constant exploration rate combined with the Boltzmann explorers.
\begin{table}[H]
  \caption{Action Exploration Heuristics}
  \centering
  \label{tab:act_expl}
  \begin{tabular}{ll|ll|ll}
    \toprule
    $\epsilon$-greedy   & $\epsilon(k)$ & Boltzmann    & $\tau(k)$    & Max-Boltz           & $(\epsilon, \tau(k))$\\
    \midrule
    $\epsilon_{Con}$    & $0.01$          & $Boltz_{Con}$& $5$           & $MB_{Con}$          & $(0.05, 5)$\\
    $\epsilon_{Lin}$    & $0.05(K-k)/K$   & $Boltz_{Lin}$& $5(K-k)/K$  & $MB_{Lin}$          & $(0.05, 5(K-k)/K$\\
    $\epsilon_{Exp}$    & $0.9995^k$      & $Boltz_{Exp}$& $5\times 0.9999^k$  & $MB_{Exp}$          & $(0.05, 5\times 0.9999^k)$\\
    \bottomrule
\end{tabular}
\end{table}

\begin{figure}[H]
  \centering
  \begin{minipage}[b]{0.49\textwidth}
    \includegraphics[width=\textwidth]{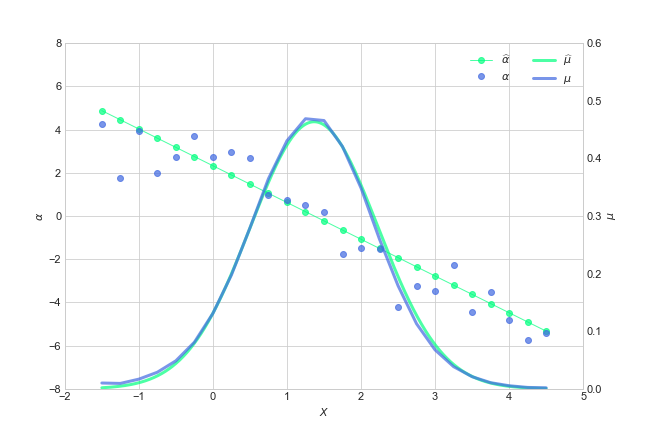}
    \caption{MFCG three-timescale Q-learning}
    \label{fig:mfcg_rates}
  \end{minipage}
  \begin{minipage}[b]{0.49\textwidth}
    \includegraphics[width=\textwidth]{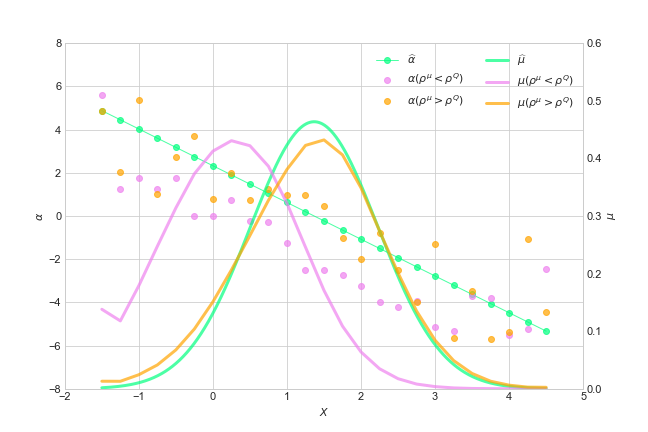}
    \caption{Two-timescale Q-learning}
    \label{fig:mfg_mfc_rates}
  \end{minipage}
\end{figure}

Algorithm \ref{alg:rl_infinite} learns the solution of the mean field control game based on three different learning rates for the Q-table and local/global distributions. Figure  \ref{fig:mfcg_rates} shows the results obtained when the learning rate parameters $(\omega_{\mu}, \omega_{Q}, \omega_{\widetilde{\mu}})$ are equal to $(0.75, 0.55, 0.15)$. The $x$-axis
represents the state variable $x$ while the left, right $y$- axes correspond to the action $\alpha(x)$ and the probability mass $\hat\mu(x)$ respectively. The green dot-marked line and continuous curves show the theoretical
solutions of the MFCG discussed in section \ref{subsec: benchmark} in terms of the control function and the asymptotic distribution at equilibrium. The blue dots and curve are the corresponding action and distribution learned by the algorithm, averaged over the last 5k episodes. Only the global distribution is plotted because the local distribution perfectly aligns with it.

Figure~\ref{fig:mfg_mfc_rates} shows how different choices of the learning rate parameters let the algorithm converge to different solutions. The green set of line and curve refers to the same theoretical solution to the MFCG problem as in Figure~\ref{fig:mfcg_rates}. The violet (resp. orange) set shown is obtained when $\omega_{\mu}=\omega_{\widetilde{\mu}}$ and their values are set to $0.75$ (resp. $0.15$) such that $\rho^{\mu} < \rho^{Q}$ (resp. $\rho^{\mu} > \rho^{Q}$). The values of actions and distributions plotted are the average of the last 5k episodes. These choices reduce the algorithm to the two-timescale approach discussed in \cite{AFL2021}. The algorithm then converges to the corresponding MFG and MFC versions of our model depending on the choice of the learning rates, where the support of the MFG deviates from the current trimmed state space.

\begin{figure}[H]
  \centering
  \begin{minipage}[b]{0.49\textwidth}
    \includegraphics[width=\textwidth]{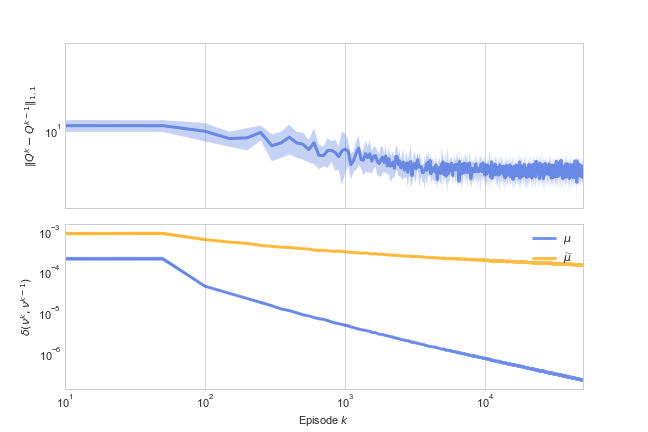}
    \caption{Total variations of $Q$, $\mu$, and $\widetilde{\mu}$}
    \label{fig: Q_norm_mu_tv}
  \end{minipage}
  \begin{minipage}[b]{0.49\textwidth}
    \includegraphics[width=\textwidth]{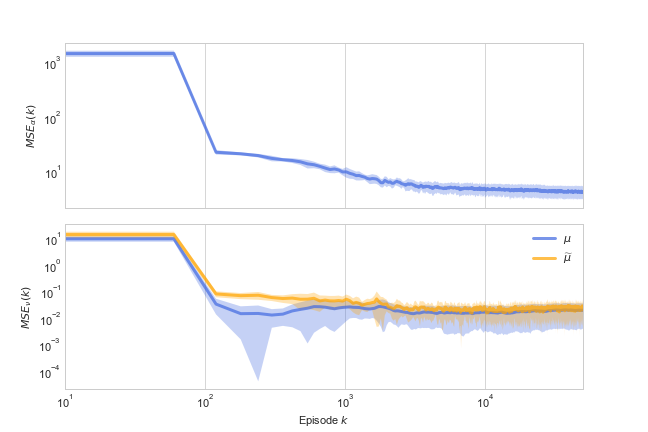}
    \caption{Mean squared errors of $\alpha$, $\mu$, and $\widetilde{\mu}$ }
    \label{fig: mse_control_distribution}
  \end{minipage}
\end{figure}

The convergence of the algorithm \ref{alg:rl_infinite} is analyzed in terms of the evolution of the estimations of the optimal $Q$ table and the local/global distributions at equilibrium w.r.t. the learning episodes. The changes are evaluated through the total variation and the $1,1$-norm as follows
\begin{align}
    \delta(\nu^{k},\nu^{k}) = \sum_{x_i \in \X} \abs{\nu^{k}(x_i)-\nu^{k-1}(x_i)}, &&
    \Vert Q^k - Q^{k-1} \Vert_{1,1} = \sum_{i,j} \abs{Q_{i,j}^k - Q_{i,j}^{k-1}},
\end{align}
where the episode is tracked by the index $k$ and $\nu$ is replaced by $\mu$ and $\widetilde{\mu}$.
Figure \ref{fig: Q_norm_mu_tv} shows how the convergence improves w.r.t. the number of episodes. The $x$-axis represents the learning episode $k$. The $y$-axis represents the value of the $1,1$-norm and the total variation respectively with the averaged values over 10 runs (solid line) and standard deviations (shaded region).

The optimal control function learned by the algorithm is evaluated w.r.t. the limiting distribution of the population at the equilibrium. In particular, we analyze the mean square error averaged over multiple runs as follows
\begin{align}
\text{MSE}_{\hat\alpha}(i,k) = \sum_{j=0}^{|\X|-1} (\alpha^{i,k}(x_j)-\hat\alpha(x_j))^2 \hat\mu(x_j),
&& \text{MSE}_{\hat\alpha}(k) = \frac{1}{\#runs}\sum_{i=0}^{\#runs} \text{MSE}_{\hat\alpha}(i,k), 
\end{align}
where $\hat\mu(x_j) = \int_{x_{j-1}}^{x_j}d\mu(x)$ is obtained by the asymptotic distribution at equilibrium $\mu$ using the convention $x_{-1}=-\infty$.
Similarly, we evaluate the learning of the first moment of the asymptotic distribution at equilibrium as 
\begin{align}
\text{MSE}_{\overline\mu}(k) = \frac{1}{\#runs}\sum_{i=0}^{\#runs} (\overline\mu_T^{i,k}- \overline{\hat{\mu}})^2,
&&\text{MSE}_{\overline{\widetilde{\mu}}}(k) = \frac{1}{\#runs}\sum_{i=0}^{\#runs} (\overline{\widetilde{\mu}}_T^{i,k}- \overline{\hat{\mu}})^2.
\end{align}

Figure \ref{fig: mse_control_distribution} shows the decrease of the errors w.r.t. the number of learning episodes. The $x$-axis corresponds to the learning episode $k$. The $y$-axis represents the errors
averaged over 10 runs (solid line) and their standard deviations (shaded region).


\begin{figure}[H]
  \centering
  \begin{minipage}[b]{0.49\textwidth}
    \includegraphics[width=\textwidth]{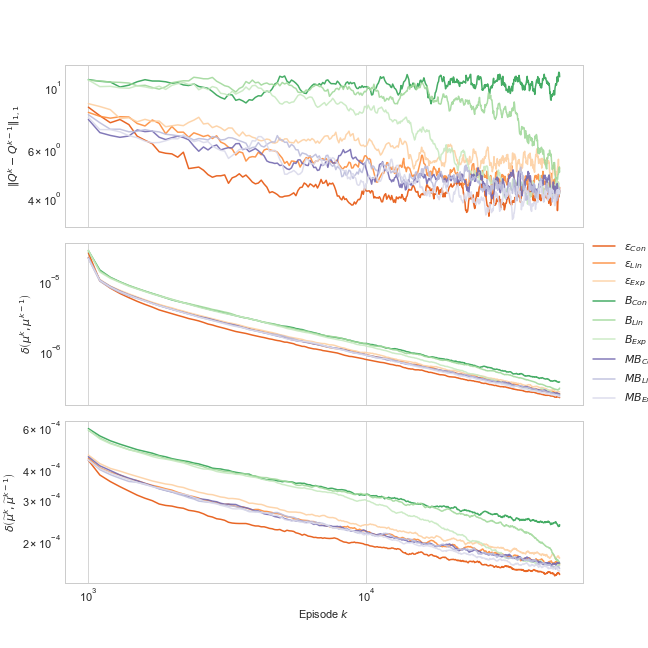}
    \caption{Comparison of total variations}
    \label{fig:act_expl_norm}
  \end{minipage}
  \begin{minipage}[b]{0.49\textwidth}
    \includegraphics[width=\textwidth]{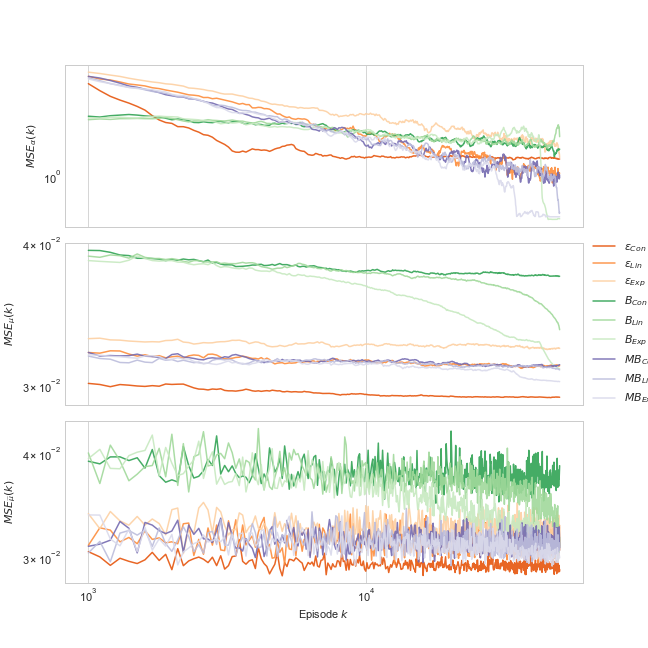}
    \caption{Comparison of mean squared errors}
    \label{fig:act_expl_mse}
  \end{minipage}
\end{figure}

We conclude by presenting an empirical comparison of the action exploration strategies discussed in section \ref{sec: action_exploration}. Figures \ref{fig:act_expl_norm} and \ref{fig:act_expl_mse} show the results obtained by applying the $\epsilon$-greedy (red set of lines), Boltzman (green set of lines), and Max-Boltzman (purple set of lines) exploration rules when the rate is constant, linear, or exponential decaying w.r.t. the episodes, as in Table~\ref{tab:act_expl}. The subplot on top of Figure~\ref{fig:act_expl_norm} shows the $1,1$-norm of the learned $Q$ table and the two subplots following are total variations of learned $\mu$ and $\widetilde{\mu}$ distribution. The $x$-axis is the log-scaled episode number, and the $y$-axes correspond to the value of those total variations. The $\epsilon$-greedy with constant exploration rate $\epsilon_{Con}$ surprisingly outperforms any other heuristic in converging speed. The worst result is obtained by the Boltzmann exploration group. $B_{Con}$ fails to converge in $Q$ table, while $B_{Lin}$ and $B_{Exp}$ waste almost 10k episodes before the variation is reduced. Obviously, the Boltzmann exploration set is under-tuned with a high initial temperature, and reducing it will hopefully improve its performance. Recall that we aim to control the exploration propensity via the probability distribution, however, the Boltzmann distribution \eqref{eq:boltzmann} depends on the value of the $Q$ table whose scale is previously unknown. Thus, it requires extra investigation to figure out both the temperature range and the decaying rate. The Max-Boltzmann exploration set performs mediocrely, which is due to its under-tuned Boltzmann component. On the contrary, in Figure~\ref{fig:act_expl_mse}, we observe that most of the heuristics result in lower mean squared error on $\alpha$ than $\epsilon_{Con}$, except $B_{Con}$, $B_{Exp}$, and $\epsilon_{Exp}$. Despite that the $\epsilon_{Con}$ still achieves the lowest mean squared error, this result indicates that the Boltzmann and Max-Boltzmann explorations could possibly lead to better results if well-tuned. Therefore, in the linear-quadratic bank borrowing and lending MFCG, the naive $\epsilon$-greedy heuristic handles the learning task well and can serve as a useful benchmark for developing more sophisticated exploration strategies.



\bibliographystyle{chicago}
\bibliography{ref}

\end{document}